\documentclass[12pt]{article}
\usepackage{amssymb,latexsym,start2e,pstcol,pst-plot,amsthm,amsmath}
\usepackage{stmaryrd}

\begin{document}
\pagestyle{plain}

\title{Counting permutations by congruence class of major index
}
\author{
H\'el\`ene Barcelo\\[-5pt]
\small  Department of Mathematics, Arizona State University.\\[-5pt]
\small Tempe, AZ 85287-1804, \texttt{barcelo@asu.edu}\\[5pt]
Bruce E. Sagan\\[-5pt]
\small Department of Mathematics, Michigan State University,\\[-5pt]
\small East Lansing, MI 48824-1027, \texttt{sagan@math.msu.edu}\\[5pt]
Sheila Sundaram\\[-5pt]
\small Department of Mathematics, Bard College,\\[-5pt]
\small Annandale-on-Hudson, NY 12504-5000, \texttt{sundaram@bard.edu}
}

\date{{\it This paper is dedicated to the memory of Bob Maule who did seminal work
    on the subject.}\\[10pt]
       \today\\[10pt]
	\begin{flushleft}
	\end{flushleft}}

\maketitle

\begin{abstract}
Consider $S_n$, the symmetric group on $n$ letters, and let
$\maj\pi$ denote the major index of $\pi\in S_n$.  Given positive
integers $k,l$ and nonnegative integers $i,j$, define
$$
m_n^{k,l}(i,j):=\#\{\pi\in S_n\ :\ \maj\pi\Cong i\ (\Mod k)
\mbox{ and } \maj\pi^{-1}\Cong j\ (\Mod l)\}
$$
We prove bijectively that if $k,l$ are relatively prime and at
most $n$ then
$$
m_n^{k,l}(i,j)=\frac{n!}{kl}
$$
which, surprisingly, does not depend on $i$ and $j$.  Equivalently, if
$m_n^{k,l}(i,j)$ is interpreted as the $(i,j)$-entry 
of a matrix $m_n^{k,l}$ then this is a constant matrix under the
stated conditions.
This bijection is extended to show the more general result  that 
for $d\ge 1$ and $k,l$ relatively prime, the matrix $m_n^{kd,ld}$ admits 
a block decompostion where each block is
the matrix $m_n^{d,d}/(kl)$.  
We also give an explicit formula for $m_n^{n,n}$ and show that if
$p$ is prime then $m_{np}^{p,p}$ has a simple block decomposition.
To prove these results, we use the representation theory of the
symmetric group and certain restricted shuffles.
\end{abstract}

\section{Introduction}

Let $S_n$ denote the symmetric group consisting of all permutations
$\pi$ of the set
$[n]:=\{1,2,\ldots,n\}$.  
If we write $\pi=a_1a_2\ldots a_n$ then the 
{\it major index of $\pi$\/} is
$$
\maj\pi:=\sum_{a_i>a_{i+1}} i.
$$

Let $k,l$ be positive integers and let $i,j$ be nonnegative integers.
We wish to study the cardinalities
$$
m_n^{k,l}(i,j) = \# \{\pi\in S_n\ :\ \maj\pi\Cong i\ (\Mod k)
\mbox{ and } \maj\pi^{-1}\Cong j\ (\Mod l)\}.
$$
We will often omit the superscript $k,l$ both for readability and because
the parameters will be clear from context.  Note that directly from
the definition we have $m_n^{k,l}(i,j)=m_n^{l,k}(j,i)$.
One of our main objectives is to give a bijective proof of the
following theorem. 
\bth
\label{gr}
Let $k,l$ be relatively prime and less than or equal to $n$.  Then
$$
m_n^{k,l}(i,j) = \frac{n!}{kl}.
$$
\eth

This theorem is striking because the right-hand side of the equality
does not depend on $i,j$.  It is easy to prove algebraically based on
results of Gordon~\cite{gor:ttm} and Roselle~\cite{ros:cae}.  This was
done in a paper of Barcelo, Maule, and Sundaram~\cite{bms:cpp} where they
also provided combinatorial proofs of special cases of this result.
Here we will give a bijective proof
with no restrictions other than those in the statement of the theorem.
These restrictions are necessary since the result is no longer true
without them.  However, we will generalize our bijection to cover the case
where the moduli are allowed to have a common factor and prove the
following.
\bth
Let $k,l$ be relatively prime and let $d\ge1$.  Then
$$
m_n^{kd,ld}(i,j)=\frac{m_n^{d,d}(i,j)}{kl}.
$$
\eth

The second half of the paper will be devoted to investigating
$m_n^{k,l}(i,j)$ when $k=l$.  We will give an explicit formula for
$m_n^{n,n}(i,j)$.  It will also be shown that when $k=l$ 
is a prime power, the matrix whose $(i,j)$th entry is $m_n^{k,k}(i,j)$ has a
nice block form.  Our tools will include restricted shuffles
and results from the
representation theory of the symmetric group.

\medskip

{\it Dedication.}
Bob Maule was an accomplished actuary who took an early retirement
due to health problems. After a few years, he enrolled in the graduate program
in mathematics at Arizona State University, for the sheer pleasure of doing
mathematics.

During a combinatorics course taught by H\'el\`ene Barcelo,
Bob became quite interested in a question which she was investigating with
Sheila Sundaram, namely, the
distribution of the values of the major index among the permutations of $S_n$.
He worked relentlessly for several months, analysing the patterns
that were slowly emerging from his computations. This is how he came
to develop a matrix approach that we further developed and use here.
His contributions to the solution of the original problem were
crucial, and led to a joint
paper with H\'el\`ene and Sheila~\cite{bms:cpp}.
Despite serious health problems
he continued his computations which foreshadowed
several of the results in Section~\ref{ppp}.
Unfortunately, Bob passed away before he could properly formalize his ideas.

It is undoubtedly his enthusiasm and hard work that brought H\'el\`ene
back to this subject. She always felt that it was an honor to be his advisor.
He was an exceptional person,  a constant source of inspiration and
a very enjoyable person to work with.
He deserves our respect both as a person and as a mathematician in the becoming.

In tribute to Bob's inspirational work we are happy to dedicate
this article to his memory, hoping that in doing so he will be remembered in
our community as a graduate student whose sole motivation was the pleasure of
doing mathematics.

\section{Preliminaries}

Before embarking on a proof of Theorem~\ref{gr}, we would like to
restate it in a form more amenable to bijective arguments.
To do so, we will also need another common combinatorial statistic.
The {\it inversion number\/} of $\pi\in S_n$ is
$$
\inv\pi:=\#\{(a_i,a_j)\ :\ i<j \mbox{ and } a_i>a_j\}.
$$
Foata and Sch\"utzenberger~\cite{fs:mii} proved bijectively that the
statistics $\maj$ and $\inv$ are equidistributed over $S_n$.  In fact,
their bijection also shows that the joint distribution of the pair
$(\maj\pi,\maj\pi^{-1})$ is the same as that of $(\inv\pi,\maj\pi)$.
Also, it is trivial to see that $\inv\pi^{-1}=\inv\pi$.  It follows
that
$$
m_n^{k,l}(i,j) = \#\{\pi\in S_n\ :\ \inv\pi\Cong i\ (\Mod k)
\mbox{ and } \maj\pi^{-1}\Cong j\ (\Mod l)\}
$$
and this is the combinatorial interpretation for these numbers that we
will use for most of the rest of the paper.  We will also need the
corresponding sets
$$
M_n^{k,l}(i,j) = \{\pi\in S_n\ :\ \inv\pi\Cong i\ (\Mod k)
\mbox{ and } \maj\pi^{-1}\Cong j\ (\Mod l)\}.
$$
It will often be convenient to think of these as
the $(i,j)$ entries of matrices $m_n^{k,l}$ and $M_n^{k,l}$,
respectively.  

To see how this change of viewpoint simplifies things, we will
give a bijective proof of a weaker form of Theorem~\ref{gr}
where we only consider one of the two statistics.  
We will also need this result in the proof of the theorem itself.
Another
combinatorial proof of this result can be found in~\cite{bms:cpp}, but
ours has the advantage of being simpler and not using induction.
Let 
\bea
m_n^{k}(i)
	&:=&\#\{\pi\in S_n\ :\ \maj\pi\Cong i\ (\Mod k)\}\\
	&=&\#\{\pi\in S_n\ :\ \maj\pi^{-1}\Cong i\ (\Mod k)\}\\
	&=&\#\{\pi\in S_n\ :\ \inv\pi\Cong i\ (\Mod k)\}.
\eea
\bpr
\label{mnk}
If $k\le n$ then
$$
m_n^{k}(i) =\frac{n!}{k}
$$
\epr
\pf\
It suffices to show that $S_n$ can be partitioned into 
subsets of the form $S=\{\pi_0,\ldots,\pi_{k-1}\}$ where
\beq
\label{pii}
\inv \pi_i\Cong \inv\pi_0 + i (\Mod k)
\eeq  
for $0\le i<k$.  
Given $\pi=a_1a_2\ldots a_n\in S_n$, we construct the subset $S$
containing it as follows. Let $a_m$ be the maximal element of the
prefix $a_1 a_2 \ldots a_k$ of $\pi$.  Let $\si$ be the sequence
formed from $\pi$ by removing $a_m$.  Finally form $\pi_i$ by
inserting $a_m$ into the $i$th space of $\si$, where the space
completely to the left of $\si$ is counted as space 0.  It is easy to
see that equation~\ree{pii} holds, so we are done.
\Qqed

\section{Proof of Theorem~\ref{gr}}

In order to prove Theorem~\ref{gr} we will need a nice combinatorial
interpretation of $\maj\pi^{-1}$ which we will henceforth write as
$\imaj\pi$.   In fact, it follows immediately from the definitions
that
$$
\imaj\pi =\sum_{\mbox{\scriptsize $i+1$ left of $i$ in $\pi$}} i
$$
and that is how the reader should think of calculating this number.

Given $n$ and $l$ with $n>l$ we will also need a particular bijection
$f=f_l$ from $S_n$ to itself defined as follows.  The reader
may wish to also read the example at the end of the paragraph while they read
the definition.   If
$\tau=a_1 a_2\ldots a_n$ then let $I=\{i_1<i_2<\ldots<i_l\}$ be the
indices such that $\pi=a_{i_1} a_{i_2}\ldots a_{i_l}$ is a permutation
of $[l]$.  Let $\si$ be the subsequence of $\pi$ indexed by $[n]-I$, thus
$\tau$ is a shuffle of $\pi$ and $\si$.  Consider
$J=\{i_1+1,i_2+1,\ldots,i_l+1\}$ where the sums are taken modulo $n$.
Then define $\tau'=f_l(\tau)$ to be the shuffle of $\pi$ and
$\si$ such that $\tau'$ restricted to $J$ and $[n]-J$ are $\pi$ and
$\si$, respectively.  Note that $f_l$ is clearly bijective since one
can construct its inverse in exactly the same manner by just
subtracting one from each element of $I$.
By way of illustration, suppose that $n=7$, $l=4$, and
$\tau=6 3 7 1 4 5 2$.  So $I=\{2,4,5,7\}$ corresponding to
$\pi=3 1 4 2$.  Also $\si= 6 7 5$.  Thus $J=\{1,3,5,6\}$ and
$\tau'=3 6 1 7 4 2 5$.
 
If $\tau$ is a shuffle of $\pi$ and $\si$ then it will be useful to
let $\inv_\tau (\pi,\si)$ denote the number of inversion pairs in $\tau$
with one element of the pair in $\pi$ and the other in $\si$.
If $\tau'=f_l(\tau)$, we
claim that
\beq
\label{invtau}
\inv_{\tau'}(\pi,\si)=
\case{\inv_\tau(\pi,\si)+l}{if $n\not\in I$,}
{\inv_\tau(\pi,\si)-(n-l)}{if $n\in I$.}
\eeq
To see this, note that since every element of $\pi$ is less than every
element of $\si$, then an inversion is created every time an element
of $\si$ preceeds an element of $\pi$.  If $n\not\in I$ then there is
no wrap-around when passing from $\tau$ to $\tau'$ and so there is one
new inversion created for each of the $l$ elements of $\pi$.  If $n\in
I$ then an element of $\pi$ is moved from the back of $\tau$ to the
front, so the inversions of this element
with the $n-l$ elements of $\si$ are lost.

We also define $\imaj_\tau(\pi,\si)$ as the subsum
of $\imaj\tau$ over those pairs $(i,i+1)$ with $i\in\pi$ and
$i+1\in\si$ or vice-versa.  Note that for the shuffles considered two
paragraphs ago,
\beq
\label{imajtau}
\mbox{$\imaj_\tau(\pi,\si)=0$ or $l$.}
\eeq

As a first step toward proving Theorem~\ref{gr}, we prove the
following special case.
\ble
\label{grbase}
Let $l$ be less than and relatively prime to $n$.  Then
$$
m_n^{n,l}(i,j) = \frac{n!}{nl}.
$$
\ele
\bprf
We claim that 
the map $f=f_l$ restricts to a bijection from
$M_n(i,j)$ to $M_n(i+l,j)$.  Keeping
the notation from the definition of $f$ and using equation~\ree{invtau},
we see that
\bea
\inv\tau'&=&\inv\pi+\inv\si+\inv_{\tau'}(\pi,\si)\\
&\Cong&\inv\pi+\inv\si+\inv_\tau(\pi,\si)+l\ (\Mod n)\\
&\Cong&\inv\tau+l\ (\Mod n).
\eea
Thus the row indices in $M_n$ change as desired.  For the columns
note that, by equation~\ree{imajtau},
\bea
\imaj\tau'&=&\imaj\pi+\imaj\si+\imaj_{\tau'}(\pi,\si)\\
&\Cong&\imaj\pi+\imaj\si\ (\Mod l)\\
&\Cong&\imaj\tau\ (\Mod l).
\eea
Hence $f$ restricts as claimed.

Since $l$ is relatively prime to $n$, the set of multiples of $l$
intersects every congruence class modulo $n$.  So iterating $f$ will establish
bijections between the sets $M_n(1,j),M_n(2,j),\ldots,M_n(n,j)$ for
any $j$.  But then by Proposition~\ref{mnk} we must have
$$
m_n^{n,l}(i,j)=\frac{m_n^{l}(j)}{n}=\frac{n!}{nl}
$$
as desired.
\eprf

The previous lemma will form the base case for an inductive proof of
Theorem~\ref{gr}.  For the induction step, we will need a restricted
type of shuffle.  We let $\pi\shu\si$ denote the set of shuffles
of the sequences $\pi$ and $\si$, e.g.,
$$
12\shu 43=
\{1243,\ 1423,\ 1432,\ 4123,\ 4132,\ 4312\}.
$$
We extend this notation (and all future variants of it) to sets by
letting $M\shu N=\cup(\pi\shu\si)$ where the union is over all 
$\pi\in M$ and $\si\in N$.

In order to get a permutation from shuffling two permutations,
define $\pi\shu^+\si=\pi\shu\tau$ where $\tau$ is the sequence formed
by adding $|\pi|$ to every element of $\si$.  For example 
$12\shu^+ 21$ would give the same set as displayed above.  

We will also need to specify where the copy of $\pi$ (and thus of
$\si$) occurs in a shuffle.  So given $M\sbe S_l$, $N\sbe S_{n-l}$,
and $I=\{i_1,i_2,\ldots,i_l\}\sbe[n]$, we define
$M\shu_I^+ N$ to be all permutations  $\tau=a_1a_1\ldots a_n\in M\shu^+ N$ such that
$a_{i_1}a_{i_2}\ldots a_{i_l}\in M$.  By way of illustration
$$
\{12,\ 21\} \shu_{\{2,5\}}^+ \{231,\ 321\}
=\{41532,\ 42531,\ 51432,\ 52431\}.
$$
Note that if $\pi\in M$ and $\si\in N$ then $\inv_\tau(\pi,\si)$ is
constant for all $\tau\in M\shu_I^+ N$.  So in this setting define the
{\it weight of $I$\/}, $\wt I$, to be this constant value
$$
\wt I = \sum_{i\in I} i - {l+1\choose 2}.
$$

The next result will permit us to complete the proof of
Theorem~\ref{gr}.  In it, $\uplus$ denotes disjoint union.
\ble
\label{grind}
Given $n,k,l$ with $n\ge l$, we have
$$
M_n(i,j)=\biguplus\left[M_l(i',j')\shu_I^+ M_{n-l}(i'',j'')\right]
$$
where the disjoint union is over all $i',j',i'',j'',I$ such that
$$
i\Cong i'+i''+\wt I\ (\Mod k)\qmq{and} j\Cong j'+j''\ (\Mod l).
$$
\ele
\bprf
To show that the right-hand side is contained in the left, let
$\tau=\pi\shu_I^+\si$ where $\pi\in M_l(i',j')$ and 
$\si\in M_{n-l}(i'',j'')$.  Then
$$
\inv\tau=\inv\pi+\inv\si+\inv_\tau(\pi,\si)=i'+i''+\wt I
\Cong i\ (\Mod k).
$$
Also, since $\pi\in S_l$, we have
$$
\imaj\tau=\imaj\pi+\imaj\si+\imaj_\tau(\pi,\si)
\Cong j'+j''\Cong j\ (\Mod l).
$$
Thus $\tau\in M_n(i,j)$.

To show the reverse containment, suppose $\tau\in M_n(i,j)$ is given.
Let $I$ be the 
indices where the elements of $[l]$ appear in $\tau$.  Also let
$\pi$ and $\si$ be $\tau$ restricted to $I$ and to $[n]-I$, respectively
(where $l$ has been subtracted from every element of the latter
restriction). Then $\tau=\pi\shu_I^+ \si$ with $\pi\in M_l(i',j')$
and $\si\in M_{n-l}(i'',j'')$ for $i',i'',j',j''$ satisfying the
equations in the statement of the Lemma.  Hence we are done.
\eprf

As an application of the two previous lemmas, let us reprove a result
from~\cite{bms:cpp} which we will need later.
\bco
\label{n+1}
We have
$$
m_{n+1}^{n,n}(i,j)=m_n^{n,n}(i,j)+(n-1)!.
$$
\eco
\bprf
In Lemma~\ref{grind}, replace $n$ by $n+1$ and let $k=l=n$.  Note
that $M_1(i'',j'')=\{1\}$ when $i''=j''=0$ and is the empty set
otherwise.  So taking cardinalities and using Lemma~\ref{grbase} gives
$$
m_{n+1}(i,j)=\sum_{i',I} m_n(i',j)
=m_n(i,j)+\sum_{i'=0}^{n-1} m_n(i',j)
=m_n(i,j)+(n-1)!
$$
as desired.
\eprf

We now have all the tools in place to prove Theorem~\ref{gr} which we
restate here for convenience.
\bth
\label{main}
Let $k,l$ be relatively prime and less than or equal to $n$.  Then
$$
m_n^{k,l}(i,j) = \frac{n!}{kl}.
$$
\eth
\bprf
We proceed by induction on $n$.  As noted after the original definition, the matrix
$m_n^{k,l}$ is the transpose of $m_n^{l,k}$.  So it will suffice to prove the
result for either $l<k$ or $l>k$ in both the base case and induction step.

We take as our base case when $l<k=n$.  But this has already been
taken care of by Lemma~\ref{grbase}.

For the induction step we assume $k<l<n$.  Since $l>k$ we can appeal
to the base case to see that
$m_l(i,j)$ is independent of $i,j$.  Thus there are bijections
between any two sets of the form $M_l(i,j)$.  We will now
construct a bijection between $M_n(i,j)$ and $M_n(i,j+1)$.
By Proposition 2.1, this will finish the proof.
Decompose these two sets as in
Lemma~\ref{grind}.  Then, since $i$ is being held constant, every set
$M_{n-l}(i'',j'')$ appearing in the expansion of the $M_n(i,j)$ also occurs
in that of $M_n(i,j+1)$.  The only difference is that in the first
expansion it is shuffled with $M_l(i',j')$ and in the second with
$M_l(i',j'+1)$.   But there is a bijection between these two sets and
so also between the corresponding shuffles.  It follows that we have a
bijection between the disjoint unions, i.e., between $M_n(i,j)$ and $M_n(i,j+1)$.
\eprf

As already mentioned, the previous theorem is not true as stated if
$k,l$ are 
not relatively prime and we shall see some examples of this in the
next section dealing with the case when $k=l$.  However, one can
extend this result to the case where the parameters have a greatest
common divisor $d$ for any $d\ge1$ as follows.
\bth
\label{dthm}
Let $k,l$ be relatively prime and let $d\ge1$.  Then
$$
m_n^{kd,ld}(i,j)=\frac{m_n^{d,d}(i,j)}{kl}.
$$
Otherwise put, $M_n^{kd,ld}$ admits a block decomposition into $kl$
submatrices of dimension $d\times d$, each of which equals
$\frac{1}{kl}M_n^{d,d}$.
\eth
\bprf
The proof is similar to that of Theorem~\ref{main}.  So we will merely
sketch it, adding details only when there are significant differences
from the previous proof.

We first need an analogue of Proposition~\ref{mnk} which is that
\beq
\label{mnkeq}
m_n^{kd,d}(i,j)=\frac{m_n^{d,d}(i,j)}{k}
\eeq
for any $k$ with $n\ge kd$.  To show this, it suffices to find a
bijection $g:M_n^{kd,d}(i,j)\ra M_n^{kd,d}(i+d,j)$ for all $i,j$ since
then
$$
m_n^{d,d}(i,j)=\sum_{i'} m_n^{kd,d}(i',j) = k\ m_n^{kd,d}(i,j)
$$
where the sum is over all $i'$ with $i'\Cong i\ (\Mod d)$ and 
$1\le i'\le kd$.  Given $\tau\in M_n^{kd,d}(i,j)$, write 
$\tau=\pi\shu_I^+ \si$ where $\pi$ is the subsequence of $\tau$ which
is a permutation of $[kd]$, which also uniquely determines $\si$ and $I$.  Define
$g(\tau)=f_d(\pi)\shu_I^+ \si$ where $f_d$ is the map defined at the
beginning of this section.  It is easy to verify, using
computations similar to those in Lemma~\ref{grbase}, that 
$g(\tau)\in M_n^{kd,d}(i+d,j)$ and that $g$ is invertible.  So this
proves equation~\ree{mnkeq}.

Next we need a version of Lemma~\ref{grbase} itself which is
\beq
\label{grbaseeq}
m_{kd}^{kd,ld}(i,j)=\frac{m_{kd}^{d,d}(i,j)}{kl}
\eeq
for $k,l$ relatively prime and $l<k$.  Recalling that $m_n^{k,l}$ and
$m_n^{l,k}$ are transposes, equation~\ree{mnkeq} can be rewritten
$$
m_n^{d,ld}(i,j)=\frac{m_n^{d,d}(i,j)}{l}.
$$
We also have
$$
m_n^{d,ld}(i,j)=\sum_{i'} m_n^{kd,ld}(i',j)
$$
where the sum is over all $i'$ with $i'\Cong i\ (\Mod d)$ and 
$1\le i'\le kd$.  So to prove equation~\ree{grbaseeq}, it suffices to
find a bijection between $M_{kd}^{kd,ld}(i,j)$ and
$M_{kd}^{kd,ld}(i+d,j)$ for all $i,j$.  Using $f_{ld}$
gives a bijection between 
$M_{kd}^{kd,ld}(i,j)$ and $M_{kd}^{kd,ld}(i+ld,j)$.  But since $k$ and
$l$ are relatively prime, iteration of this map eventually produces
the desired bijection and we have equation~\ree{grbaseeq}.

Finally, we need an induction on $n$ to prove the full result, where
the previous paragraph gives us the base case when $n=kd$ (assuming,
without loss of generality, that $l<k$).  But now Lemma~\ref{grind}
can be used in much the same way as in the proof of
Theorem~\ref{main} to complete the induction step.  Specifically, this
Lemma can be used to lift the bijections in $M_{kd}^{kd,ld}$ to
bijections between $M_n^{kd,ld}(i,j)$ and $M_n^{kd,ld}(i,j+d)$ for all
$i,j$.  By transposition, we also get bijections between
$M_n^{kd,ld}(i,j)$ and $M_n^{kd,ld}(i+d,j)$.  Then the proof is
finished by noting that
$$
m_n^{d,d}(i,j)=\sum_{i',j'} m_n^{kd,ld}(i',j')
$$
where the sum is over all $i',j'$ with $i'\Cong i\ (\Mod d)$, 
$j'\Cong j\ (\Mod d)$, $1\le i'\le kd$, and $1\le j'\le ld$.
\eprf

\section{The case $k=l$}

Theorem~\ref{dthm} reduces computation of the matrices $m_n^{k,l}$ to
the case where $k=l$ which we will now consider.  We first derive a
formula for the special case $m_n^{n,n}(i,j)$.  We will use techniques
from the representation theory of the symmetric group $S_n$.
More information about these methods can be found in the texts of
Sagan~\cite{sag:sym} or Stanley~\cite[Chapter 7]{sta:ec2}.  To state
our result, we use $\mu$ and $\phi$ for the number-theoretic
M\"obius and Euler totient functions, respectively.  We also let
$i\mt j$ denote the greatest common divisor of $i$ and $j$.

\bth
\label{base}
Let $1\le i,j\le n$.  Then
\beq
\label{mnnn}
m_n^{n,n}(i,j)=\frac{1}{n^2}\sum_{d|n} d^{n/d}(n/d)!\phi(d)^2
\frac{\mu\left(\frac{d}{i\mt d}\right)\mu\left(\frac{d}{j\mt d}\right)}
{\phi\left(\frac{d}{i\mt d}\right)\phi\left(\frac{d}{j\mt d}\right)}
\eeq
\eth
\bprf
Let $\om$ be a primitive $n$th root of unity and consider the character
$\om^i$ of the cyclic subgroup of $S_n$ generated by an $n$-cycle.
Let $\chi_{n,i}$ denote the character obtained by inducing $\om^i$ up
to $S_n$.  It is easy to see that $\chi_{n,i}$ is only nonzero on
conjugacy classes of type $d^{n/d}$ where $d|n$.  On these classes,
Foulkes~\cite{fou:csg} showed that its value is
$$
\frac{1}{n} d^{n/d}(n/d)!\phi(d)
\frac{\mu\left(\frac{d}{i\mt d}\right)}{\phi\left(\frac{d}{i\mt d}\right)}.
$$
It follows that the inner product $\ipr{\chi_{n,i}}{\chi_{n,j}}$ is
given by the right-hand side of equation~\ree{mnnn}.

The following fact was discovered independently by Kra\'skiewicz and
Weyman~\cite{kw:aca}, and by Stanley~\cite[Exercise 7.88 b]{sta:ec2}.
The multiplicity in $\chi_{n,i}$ of the 
irreducible character of $S_n$ indexed by a partition $\la$ is
the number, $f^\la_{n,i}$, of standard Young tableaux of shape $\la$
with major index congruent to $i$ modulo $n$.  Using the
decomposition into irreducibles we obtain
$$
\ipr{\chi_{n,i}}{\chi_{n,j}}=\sum_{\la\ptn n}f^\la_{n,i} f^\la_{n,j}.
$$
But via the Robinson-Schensted correspondence one sees that this sum is
exactly $m_n^{n,n}(i,j)$, so we are done.
\eprf

As an application of this theorem, we note the following useful
result.
\bco
\label{gcd}
Suppose $i\mt n=i'\mt n$ and $j\mt n=j'\mt n$ where $1\le i,i',j,j'\le
n$.  Then
$$
m_n^{n,n}(i,j)=m_n^{n,n}(i',j').
$$
So to determine the matrix $m_n^{n,n}$ it suffices to determine the
entries $m_n^{n,n}(i,j)$ when $i$ and $j$ divide $n$.  Similarly, the
numbers $m_{n+1}^{n,n}(i,j)$ only depend on $i\mt n$, $j\mt n$, and $n$.
\eco
\bprf
The first part follows immediately from equation~\ree{mnnn} and the
fact that we have $i\mt d=(i\mt n)\mt d$ whenever $d|n$.  The second part
follows from the first and Corollary~\ref{n+1}.
\eprf

The previous theorem will be useful for the base cases of various
inductive proofs.  To do the induction step, we will need a type of
shuffle $\pi\shu_\ga^+ \si$ where the distances between elements of
$\pi$ are restricted.  Suppose $\ga=(g_1,g_2,\ldots,g_{l-1})$ is a composition
(ordered partition) and $l=|\pi|$.  Let
$\pi\shu_\ga^+ \si$ be the set of all shuffles in $\pi\shu^+\si$ such
that if the elements of $\pi$ are at indices $I=\{i_1<i_2<\ldots<i_l\}$
then $i_{t+1}-i_t=g_t$ for $1\le t< l$.  By way of illustration,
$$
132\shu_{(1,2)}^+ 321 =\{136254,\ 613524,\ 651342\}.
$$
In this situation, we also define the {\it weight\/} of the composition to be
$$
\wt\ga=\inv_\tau(\pi,\si)=\sum_{t=1}^{l-1} (g_t-1)(l-t)
$$
where $\tau$ is the shuffle in $\pi\shu_\ga^+\si$ whose first element
coincides with the first element of $\pi$.  In the example above
$$
\wt(1,2)=\inv_{136254}(132,654)=1.
$$

\ble
\label{ind}
Given $n,l$ with $n\ge l$, we have
\beq
\label{indeq}
M_n^{l,l}(i,j)=\biguplus\left[M_l^{l,l}(i',j')\shu_\ga^+ M_{n-l}^{l,l}(i'',j'')\right]
\eeq
where the disjoint union is over all $i',j',i'',j'',\ga$ such that
$$
i\Cong i'+i''+\wt \ga\ (\Mod l)\qmq{and} j\Cong j'+j''\ (\Mod l).
$$
\ele
\bprf
The proof is much the same as that of Lemma~\ref{grind}, just noting
that since $|\pi|=l$ we have $\inv_\tau(\pi,\si)\Cong\wt \ga\ (\Mod l)$ for all 
$\tau\in \pi\shu_\ga^+ \si$.
\eprf

\section{Primes and prime powers}

\label{ppp}

We  now specialize to the case where $n$ is a prime $p$. Then
the sum in~\ree{mnnn} simplifies greatly and we can readily write down
the entries of the matrix.  To do so conveniently in block form, let $J_{k,l}$ be
the $k\times l$ all ones matrix.  We can also use Corollary~\ref{n+1}
to give the entries for the associated matrix from $S_{p+1}$.
\bpr
Let $p$ be prime.  Then we have 
$$
M_p^{p,p}=\frac{(p-1)!}{p}J_{p,p}+\frac{1}{p}\left[\barr{cc}
  (p-1)^2J_{1,1}   &  -(p-1)J_{1,p-1}\\[10pt]
  -(p-1)J_{p-1,1}   & J_{p-1,p-1}
\earr\right]
$$
and\\
\eqed{
M_{p+1}^{p,p}=M_p^{p,p}+(p-1)! J_{p,p}.
}
\epr

In fact, we have the same block decomposition for $M_n^{p,p}$ whenever
$n$ is a multiple of $p$ or one more than a multiple of $p$.
\bth
\label{prime}
For each prime $p$ there are nonnegtive integer sequences
$(q_n)_{n\ge1}$, $(r_n)_{n\ge1}$, and $(s_n)_{n\ge1}$ such that
$$
M_{np}^{p,p}=\left[\barr{cc}
  q_n J_{1,1}   &   r_n J_{1,p-1}\\[10pt]
  r_n J_{p-1,1}   & s_n J_{p-1,p-1}
\earr\right].
$$
The matrices $M_{np+1}^{p,p}$ have the same block decomposition (for
three other sequences).
\eth
\bprf
We will prove the result for $S_{np}$ as the one for $S_{np+1}$ is
proved similarly.  We proceed by induction on $n$, with the previous proposition
providing the base case.

Since $M_{np}$ is symmetric, the statement in the theorem is
equivalent to showing that for any $i\ge0$
and any $1\le j<p-1$ we have $m_{np}(i,j)=m_{np}(i,j+1)$.  Decompose
both sides of this equation using Lemma~\ref{ind} with $l=p$ and $n$
replaced by $np$.  Note that by the choice of $l$ and $n$, induction
applies to the matrices on the right-hand side of~\ree{indeq}.
Furthermore, since $i$ is being held constant, the same $i'$, $i''$,
and $\ga$ will occur in the expansions for the $(i,j)$ and $(i,j+1)$
entries.  Also, the equation $j\Cong j'+j''\ (\Mod p)$ has precisely
$p+1$ solutions of which exactly two involve zero because of the range
of $j$.  But the same can be said of $j+1$.  Thus when one takes
cardinalities in~\ree{indeq}, the two summations will have identical
summands.  Hence we are done.
\eprf

When $n$ is a prime power $p^r$, then the sum~\ree{mnnn} also
simplifies.  Using Corollary~\ref{gcd} and symmetry,
one sees that to determine the matrix for this value of $n$ we need
only compute $m_n^{n,n}(p^i,p^j)$ for $i\le j$.  In fact,
because of zero terms in the sum, for given $i$ all
the values for $j>i$ are equal.
\bpr
Suppose $n=p^r$ where $p$ is prime and $0\le i\le j\le r$.  The matrix
$m_n^{n,n}$ is completely determined by the values
$$
m_n^{n,n}(p^i,p^j)=\frac{1}{p^{2r}}
\sum_{k=0}^{i+1} \left( p^k\right)^{p^{r-k}}\left(p^{r-k}\right)!\phi(p^k)^2\psi(i,j,k)
$$
where\\
$\textcolor{white}{\qed}\hfill\dil
\psi(i,j,k)=\left\{
\barr{ll}
1                       &\mbox{if $k\le i$,}\\[10pt]
\dil\frac{1}{(p-1)^2}   &\mbox{if $k= i+1$ and $i=j$, }\\[20pt]
\dil\frac{-1}{p-1}      &\mbox{if $k= i+1$ and $i<j$. }
\earr
\right.
\hfill\raisebox{-35pt}{\qed}$
\epr

Similarly, Theorem~\ref{prime} and its proof can be extended without
difficulty to the following result.
\bth
Suppose $k=l=p^r$ where $p$ is prime and $0\le i\le j\le r$.  The
matrix $m_{np^r}^{p^r,p^r}$ is
completely determined by its entries in
positions $(p^i,p^j)$.  Furthermore, given $i$ all these entries for
$j>i$ are equal.  The same is true for the matrix  $m_{np^r+1}^{p^r,p^r}$.\hfill\qed
\eth

In general, it is not easy to explicitly compute the sequences
$(q_n)_{n\ge1}$, $(r_n)_{n\ge1}$, and $(s_n)_{n\ge1}$ in
Theorem~\ref{prime} because the expression~\ree{indeq} of Lemma~\ref{ind} is so
complicated.  But when $p=2$ things simplify greatly and we also get
equality of the diagonal elements which is not true in general.  

To state our results compactly, let $B_n=M_n^{2,2}$ and similarly for
the matrix $b_n$.  Also, if $\Ga$ is a set of compositions of length
$|\pi|-1$ then let 
$$
\pi\shu_\Ga^+ \si=\biguplus_{\ga\in\Ga} \left(\pi\shu_\ga^+\si\right).
$$
Finally, let 
$$
O=\{1,3,5,\ldots\}\qmq{and} E=\{2,4,6,\ldots\}.
$$

The next theorem follows easily from Lemma~\ref{ind} and induction, so
we will merely state it.  None of the sequences mentioned in this
result have been previously submitted to Sloane's
encyclopedia~\cite{slo:oei}.
\bth
Suppose $n\ge2$.
\ben
\item \raisebox{-13pt}{$\barr{rcl}
  B_n(i,j)&=&\left(12\shu_O^+ B_{n-2}(i,j)\right)\biguplus
                 \left(12\shu_E^+ B_{n-2}(i+1,j)\right)\biguplus\\[10pt]
		&& \left(21\shu_E^+ B_{n-2}(i,j+1)\right)\biguplus
		 \left(21\shu_O^+ B_{n-2}(i+1,j+1)\right).
		 \earr$}
\item \raisebox{-13pt}{$\barr{rcl}
  b_{2n}(i,j)&=&2n^2 b_{2n-2}(i,j)+2n(n-1)b_{2n-2}(i+1,j),\\[10pt]
  b_{2n+1}(i,j)&=&2n(n+1) b_{2n-1}(i,j)+2n^2 b_{2n-1}(i+1,j).
		 \earr$}
\item  The matrices $c_{2n}:=b_{2n}/(2^{n-1}n!)$ and
$c_{2n+1}:=b_{2n+1}/(2^{n-1}n!)$ have integer entries satisfying
\bea
  c_{2n}(i,j)&=&n c_{2n-2}(i,j)+(n-1)c_{2n-2}(i+1,j),\\
  c_{2n+1}(i,j)&=&(n+1) c_{2n-1}(i,j)+n c_{2n-1}(i+1,j).
\eea
\item $b_n(i,j)=b_n(i+1,j+1).$
\item $b_{2n}=(2n)b_{2n-1}$.\hfill\qed
\een 
\eth

\bigskip
\bibliographystyle{plain}
\begin{small}
\bibliography{ref}
\end{small}

\end{document}